\documentstyle[amsfonts,amssymb,12pt]{article}
\font\sixbb=msbm6
\font\eightbb=msbm8
\font\twelvebb=msbm10 scaled 1095
\newfam\bbfam
\textfont\bbfam=\twelvebb \scriptfont\bbfam=\eightbb
                           \scriptscriptfont\bbfam=\sixbb
\def\bb{\fam\bbfam\twelvebb}

\newcommand{\Int}{{\bb Z}}

\newtheorem{theorem}{\bf Theorem}[section]
\newtheorem{cclaim}[theorem]{\bf Claim}

\newtheorem{proposition}[theorem]{\bf Proposition}

\newcommand{\enp}{\begin{flushright} $\Box$ \end{flushright}}
\newcommand{\beq}[0]{\begin{equation}}
\newcommand{\enq}[0]{\end{equation}}
\newcommand{\fht}{{\widehat{f}}}

\newcommand{\zp}{\Int_p}

\newcommand{\spp}{{\rm supp}}

\newcommand{\Ght}{\widehat{G}}
\newcommand{\Hht}{\widehat{H}}

\newcommand{\HP}{H^{\perp}}
\newcommand{\GHht}{\widehat{G/H}}
\newcommand{\teta}{\tilde{\eta}}
\newcommand{\ybar}{\overline{y}}

\title{An Uncertainty Inequality for  \\ Finite Abelian Groups}
\author{ Roy Meshulam\thanks{Department of
Mathematics, Technion, Haifa 32000, Israel. e-mail:
meshulam@math.technion.ac.il}}
\begin{document}
\maketitle

\begin{abstract}
Let $G$ be a finite abelian group of order $n$. For a  complex
valued function $f$ on $G$ let $\fht$ denote the Fourier transform
of $f$. The classical uncertainty inequality asserts that if $f
\neq 0$ then
\begin{equation}
\label{clas} |\spp(f)| \cdot |\spp(\fht)| \geq |G|~~.
\end{equation}
Answering a question of Terence Tao, the following improvement of
(\ref{clas}) is shown:
\\ $~$ \\
\noindent {\bf Theorem:} Let $d_1<d_2$ be two consecutive divisors
of $n$. If $d_1 \leq k=|\spp(f)| \leq d_2$ then
\begin{eqnarray*}
|\spp(\fht)| \geq \frac{n}{d_1d_2}(d_1+d_2-k)~~.
\end{eqnarray*}
\end{abstract}

\section{Introduction}
Let $G$ be a finite abelian  group of order $n$ and let $\Ght$ be
its character group. Let $L(G)$ denote the space of complex valued
functions on $G$. For $f \in L(G)$ let  $\fht \in L(\Ght)$ denote
its Fourier transform: \begin{eqnarray*} \fht(\chi)=\sum_{x \in G}
f(x) \chi(-x)~.\end{eqnarray*} Let $\spp(f)=\{x \in G:f(x) \neq
0\}$ denote the support of $f$. The classical uncertainty
inequality (see e.g. \cite{DS89,Sm90,M92,Terras}) asserts that if
$0 \neq f \in L(G)$ then
\begin{equation}
\label{class} |\spp(f)| \cdot |\spp(\fht)| \geq n~~.
\end{equation}
For a subgroup $H<G$ let $H^{\perp}=\{ \lambda \in \Ght~:~ \ker
\lambda \supset H\}$. If $f=1_H$ is the indicator function of $H$,
then $\fht=|H|\cdot 1_{\HP}$ and (\ref{class}) is satisfied with
equality. Conversely, it can be shown (see \cite{M92}) that if $f
\in L(G)$ satisfies (\ref{class}) with equality and $f(0)=1$ then
$f(x)=1_H(x) \chi(x)$ for some $H<G$ and $\chi \in \Hht$.
\\
Recently Tao \cite{Tao} showed that (\ref{class}) can be
substantially improved when $G=\zp$  is the cyclic group of prime
order $p$.
\begin{theorem}{\cite{Tao}}
\label{modp} If $0 \neq f \in L(\zp)$ then
\begin{eqnarray*}
|\spp(f)| + |\spp(\fht)| \geq p+1~~.
\end{eqnarray*}
\end{theorem}
Tao further conjectured that one could similarly improve
(\ref{class}) for all finite abelian groups provided that
$|\spp(f)|$ stays away from any divisor of $|G|$.
\\
In this note we extend Theorem \ref{modp} to general finite abelian
groups. For an integer $n$ and a real number $1 \leq k \leq n$ let
$d_1(n,k)$ denote the largest divisor $d_1$ of $n$ such that $d_1
\leq k$, and let $d_2(n,k)$ denote the smallest divisor $d_2$ of
$n$ such that $d_2 \geq k$.
\begin{theorem}
\label{uncer} Let $f \in L(G)$ such that $0 \neq |\spp(f)|=k$ and
let $d_i=d_i(n,k)$. Then
\begin{equation}
\label{main} |\spp(\fht)| \geq \frac{n}{d_1d_2}(d_1+d_2-k)~~.
\end{equation}
\end{theorem}
{\bf Remark:} Tao noted that Theorem \ref{uncer} can also be
formulated as follows: If $f$ is a non-zero function on $G$, then
the lattice point $(|supp(f)|, |supp(\fht)|)$ lies on or above the
convex hull of the points $(|H|, |G/H|)$, where $H$ ranges over
all subgroups of $G$. The classical uncertainty inequality,
meanwhile, merely states that this lattice point lies above the
hyperbola connecting those points.
\ \\ \\
The proof of Theorem \ref{uncer} depends on Theorem \ref{modp} and
on the following two simple observations. For $~1 \leq k \leq
n=|G|~$ let
\begin{eqnarray*}\theta(G,k)=\min~\{~|\spp(\fht)|~:~0 \neq f \in
L(G)~,~|\spp(f)| \leq k~\}~~.\end{eqnarray*}
\begin{proposition}
\label{sub} Let $H$ be a subgroup of $G$ and let $1 \leq k \leq
n$. Then there exist $1 \leq s \leq |H|~$ and $~1 \leq t \leq |G/H|~$
such that $st \leq k$ and
\begin{equation}
\label{sube}
\theta(G,k) \geq
\theta(H,s)\cdot \theta(G/H,t)~~.
\end{equation}
\end{proposition}
For $1 \leq k \leq n$ let $d_i=d_i(n,k)$ and let
$u(n,k) = \frac{n}{d_1d_2}(d_1+d_2-k)~~.$
\begin{proposition}
\label{subm}
For any divisor $d$ of $n$ and for any $1 \leq s \leq d~,~1 \leq
t \leq \frac{n}{d}$
\begin{equation}
\label{subme}
u(d,s)\cdot u(\frac{n}{d},t) \geq u(n,st)~~.
\end{equation}
\end{proposition}
{\bf Proof of Theorem \ref{uncer}:} We show by induction on
$|G|=n$ that $\theta(G,k) \geq u(n,k)$ for all $1 \leq k \leq n$.
For prime $n$, this reduces to Tao's result. Otherwise let $d$ be
a non-trivial divisor of $n$ and let $H$ be a subgroup of $G$ of
order $d$.  By Proposition \ref{sub} there exist $1 \leq s \leq d$
and $1 \leq t \leq \min\{\frac{k}{s},\frac{n}{d}\}$ such that
(\ref{sube}) holds. Combining the induction hypothesis with
(\ref{subme}) and the monotonicity of $u$, we obtain
\begin{eqnarray*}
\theta(G,k) \geq \theta(H,s) \cdot \theta(G/H,t) \geq
u(d,s)\cdot u(\frac{n}{d},t) \geq u(n,st) \geq u(n,k)~~.
\end{eqnarray*}
{\enp} The proofs of Propositions \ref{sub} and \ref{subm} are
given in Sections \ref{s:sub} and \ref{s:subm}. In Section
\ref{s:nonab} we remark on a possible extension to non-abelian
groups.

\section{Subgroups and Factor Groups}
\label{s:sub} For a subgroup $H<G$ let $q:\Ght \rightarrow \Hht$
denote the restriction homomorphism.  For each $\eta \in \Hht$
choose an arbitrary but fixed $\teta \in \Ght$ such that
$q(\teta)=\eta$. Clearly $q^{-1}(\eta)=\{\teta\cdot
\lambda~:~\lambda \in \HP\}$ and $\Ght=\{\teta\cdot \lambda~:~
\eta \in \Hht~,~\lambda \in \HP\}$. For $f \in L(G)$ and $y \in G$
let $f_y \in L(H)$ be given by $f_y(z)=f(z+y)$ for all $z \in H$.
Let $\ybar=y+H$ denote the image of $y \in G$ in $G/H$. For $\eta
\in \Hht$ let $F_{\eta} \in L(G/H)$ be defined by
$F_{\eta}(\ybar)=\widehat{f_y}(\eta)\teta(-y)~~.$
It can be checked that the right-hand side indeed depends only on
$\ybar$. A character $\lambda \in \HP$ gives rise to a character
$\lambda' \in \GHht$ given by $\lambda'(\ybar)=\lambda(y)$. The
map $\lambda \rightarrow \lambda'$ is an isomorphism between $\HP$
and $\GHht$.
\begin{cclaim}
\label{factor} For $\eta \in \Hht$ and $\lambda \in \HP$
\begin{eqnarray*}\fht(\teta\cdot \lambda)=\widehat{F_{\eta}}(\lambda')~~.\end{eqnarray*}
\end{cclaim}
{\bf Proof:} Let $(G:H)=m$ and let $G=\bigcup_{i=1}^m (y_i+H)$ be
the coset decomposition of $G$. Then
\begin{eqnarray*}
\fht(\teta\cdot \lambda)=\sum_{x \in G} f(x) \teta(-x) \lambda(-x)
=\sum_{i=1}^m \sum_{z \in H} f(z+y_i)
\teta(-z-y_i)\lambda(-z-y_i)=
\end{eqnarray*}
\begin{eqnarray*}\sum_{i=1}^m (~\sum_{z \in H} f_{y_i}(z)
\eta(-z)~) \teta(-y_i)\lambda(-y_i)=\sum_{i=1}^m
\widehat{f_{y_i}}(\eta)\teta(-y_i)\lambda(-y_i)=
\end{eqnarray*}
\begin{eqnarray*}
\sum_{i=1}^m F_{\eta}(\overline{y_i})\lambda'(-\overline{y_i})=
\widehat{F_{\eta}}(\lambda')~~.
\end{eqnarray*} {\enp}
{\bf Proof
of Proposition \ref{sub}:} Let $f \in L(G)$ with $|\spp(f)|=k >0$.
Keeping the notation of Claim \ref{factor} let
\begin{eqnarray*}I=\{1 \leq i\leq m~:~\spp(f) \bigcap (y_i+H) \neq
\emptyset~\}~\end{eqnarray*} and denote $|I|=t$. Let $\eta$ be any
element of $\Hht$. If $j \not\in I$ then $f_{y_j}=0$ and therefore
$F_{\eta}(\overline{y_j})=0$. It follows that if $F_{\eta} \neq 0$
then
\begin{equation}
\label{bnd1} |\spp(\widehat{F_{\eta}})| \geq
\theta(G/H,t)~~.
\end{equation}
By averaging there exists an $i \in I$ such that
$0<|\spp(f_{y_i})| \leq \frac{k}{t}=s~$. Let
$A=\spp(\widehat{f_{y_i}}) \subset \Hht$ then
\begin{equation}
\label{bnd2}|A| \geq \theta(H,s)~.
\end{equation}
Furthermore $F_{\eta}(\overline{y_i})=\widehat{f_{y_i}}(\eta)
\cdot \teta(y_i) \neq 0$ for all $\eta \in A$. Combining Claim
\ref{factor} with (\ref{bnd1}) and (\ref{bnd2}) it follows that
\begin{eqnarray*} |\spp(\fht)|=\sum_{\eta \in
\Hht}|\spp(\widehat{F_\eta})| \geq\sum_{\eta \in A}
|\spp(\widehat{F_\eta})| \geq \theta(H,s) \cdot
\theta(G/H,t)~.\end{eqnarray*} {\enp}

\section{A Submultiplicativity Property of $u(n,k)$}
\label{s:subm}
Let $k=st$.
For $i=1,2$ let $d_i(d,s)=a_i~,~d_i(\frac{n}{d},t)=b_i~$
and  $~d_i(n,k)=c_i~$. Then
\begin{equation}
\label{mm} m_1=\max\{a_1,\frac{k}{b_2}\} \leq s \leq
\min\{a_2,\frac{k}{b_1}\}=m_2~~.
\end{equation}
We have to show that
\begin{equation}
\label{maineq} \frac{(a_1+a_2-s)(b_1+b_2-\frac{k}{s})}{a_1 a_2 b_1
b_2} \geq \frac{c_1+c_2-k} {c_1 c_2}~~.
\end{equation}
Without loss of generality we may assume $a_1 b_1 \leq a_1 b_2
\leq a_2 b_1 \leq a_2 b_2~.$
Consider three cases:
\ \\ \\
(1) $~~a_1 b_1 \leq k \leq a_1 b_2~~.$
Since both $a_1b_1$ and $a_1b_2$ are divisors of $n$ it follows
that  $a_1 b_1 \leq c_1 \leq k \leq c_2 \leq a_1 b_2~$. By
convexity it therefore suffices to show \begin{eqnarray*}
\frac{(a_1+a_2-s)(b_1+b_2-\frac{k}{s})}{a_1 a_2 b_1 b_2} \geq
\frac{a_1b_1+a_1b_2-k}{(a_1b_1)(a_1b_2)}~~ \end{eqnarray*} or
equivalently
\begin{equation}
\label{casea} a_1(a_1+a_2-s)(b_1+b_2-\frac{k}{s}) \geq
a_1a_2b_1+a_1a_2b_2-a_2 k~~.
\end{equation}
By (\ref{mm}), $a_1 =m_1 \leq s \leq m_2 = \frac{k}{b_1}~.$ By
convexity we just have to check (\ref{casea}) for the two extreme
values of $s$:
\\
(i) $s=a_1$. Then (\ref{casea}) holds with equality.
\\
(ii) $s=\frac{k}{b_1}$. Then (\ref{casea}) is equivalent to
$(k-a_1b_1) (a_2b_1-a_1b_2) \geq 0$ which clearly holds. \ \\
\\ (2) $~~a_1 b_2 \leq k \leq a_2 b_1~~.$ Arguing as in case (1)
it suffices to show \begin{eqnarray*}
\frac{(a_1+a_2-s)(b_1+b_2-\frac{k}{s})}{a_1a_2b_1b_2} \geq
\frac{a_1b_2+a_2b_1-k}{ (a_1b_2)(a_2b_1)} \end{eqnarray*} or
equivalently
\begin{equation}
\label{caseb} (a_1+a_2-s)(b_1+b_2-\frac{k}{s}) \geq a_1b_2
+a_2b_1-k~~.
\end{equation}
We check (\ref{caseb}) for
\\
(i) $s=m_1=\frac{k}{b_2}$ . Then (\ref{caseb}) is equivalent to
$(b_2-b_1)(k-a_1b_2) \geq 0~$.
\\
(ii)  $s=m_2=\frac{k}{b_1}$ . Then (\ref{caseb}) is equivalent to
$(b_2-b_1)(a_2b_1-k) \geq 0~$. \ \\ \\ (3)  $~~a_2 b_1 \leq k \leq
a_2 b_2~~.$ It suffices to show \begin{eqnarray*}
\frac{(a_1+a_2-s)(b_1+b_2-\frac{k}{s})}{a_1a_2b_1b_2} \geq
\frac{a_2b_1+a_2b_2-k}{ (a_2b_1)(a_2b_2)} \end{eqnarray*} or
equivalently
\begin{equation}
\label{casec} a_2(a_1+a_2-s)(b_1+b_2-\frac{k}{s}) \geq a_1 a_2 b_1
+a_1 a_2 b_2-a_1 k~~.
\end{equation}
We check (\ref{casec}) for
\\
(i) $s=m_1=\frac{k}{b_2}$ . Then (\ref{casec}) is equivalent to
$(a_2 b_1-a_1 b_2)(a_2 b_2 -k ) \geq 0~$.
\\
(ii)  $s=m_2=a_2$ . Then (\ref{casec}) is in fact an
equality. {\enp}
\section{Concluding Remarks}
\label{s:nonab}
We have shown that if $0<|\spp(f)|=k$ lies between two consecutive divisors
$d_1<d_2$ of  $|G|=n$ then
$|\spp(\fht)|$ is at least the weighted average $\frac{n}{d_1 d_2}(d_1+d_2-k)$
of $\frac{n}{d_1}$ and $\frac{n}{d_2}$. It would be interesting to
obtain a similar result in the following non-abelian setting:
Let $G$ be any finite group and let
$\rho_1,\ldots,\rho_t$ be the complex irreducible representations
of $G$, where $\rho_i:G \rightarrow {\rm GL}(V_i)$ and $V_i$ is a
complex vector space of dimension $d_i$. The Fourier Transform
$\fht(\rho)$ of a function $f \in L(G)$ at a representation
$\rho:G \rightarrow {\rm GL}(V)$ is given by
\begin{eqnarray*}
\fht(\rho)=\sum_{x \in G}
f(x) \rho(x^{-1}) \in {\rm End}(V)~~.\end{eqnarray*}
Let $\mu(f)=\sum_{i=1}^t d_i \cdot {\rm rank}~\fht(\rho_i)~$. The
following non-abelian extension of (\ref{class}) was noted in
\cite{M92}.
\begin{theorem}{\cite{M92}}
\label{nonab}
 For any $0 \neq f \in L(G)$
\begin{equation}
\label{non} |\spp(f)|\cdot \mu(f) \geq |G|~.
\end{equation}
\end{theorem}
It seems likely that as in the abelian case, (\ref{non}) could be
improved when $|\spp(f)|$ is far from an order of any subgroup of
$G$.
\ \\ \\
{\bf Acknowledgement:} I would like to thank Terry Tao for helpful
comments.

\end{document}